\documentclass[english]{article}
\usepackage[T1]{fontenc}
\usepackage[latin9]{inputenc}
\usepackage{geometry}
\usepackage{verbatim}
\usepackage{amsmath}
\usepackage{amssymb}
\usepackage{esint}

\makeatletter
\newcommand{\lyxaddress}[1]{
\par {\raggedright #1
\vspace{1.4em}
\noindent\par}
}

\usepackage{babel}

\usepackage{babel}

\usepackage{babel}

\makeatother

\usepackage{babel}

\begin{document}

\title{Summed series involving $\,_{1}F_{2}$ hypergeometric functions}

\author{Jack C. Straton}

\maketitle
%\part{}

\lyxaddress{Department of Physics, Portland State University, Portland, OR 97207-0751,
USA}

\lyxaddress{straton@pdx.edu}
\begin{abstract}
In a prior paper we found that the Fourier-Legendre series of a Bessel
function of the first kind $J_{N}\left(kx\right)$ and of a modified
Bessel functions of the first kind $I_{N}\left(kx\right)$ lead to
an infinite set of series involving $\,_{1}F_{2}$ hypergeometric
functions (extracted therefrom) that could be summed, having values
that are inverse powers of the eight primes $1/\left(2^{i}3^{j}5^{k}7^{l}11^{m}13^{n}17^{o}19^{p}\right)$
multiplying powers of the coefficient $k$, for the first 22 terms in each series. The present paper shows
how to generate additional, doubly infinite summed series involving $\,_{1}F_{2}$ hypergeometric
functions from Chebyshev polynomial expansions of Bessel functions,
and trebly infinite sets of summed series involving $\,_{1}F_{2}$
hypergeometric functions from Gegenbauer polynomial expansions of
Bessel functions.
\end{abstract}
%\pacs{00.00, 0.00, 0.0}
\vspace{2pc}

\section{Introduction}

In a prior paper \cite{stra24b} we refined Keating's \cite{KeatingPhD}
derivation of the coefficient set of the Fourier-Legendre series for
the Bessel function $J_{N}\left(kx\right)$ to be

\begin{equation}
J_{N}\left(kx\right)=\sum_{L=0}^{\infty}a_{LN}\left(k\right)P_{L}(x)\quad\label{eq:Fourier-Legendre series}\end{equation}
where

\begin{equation}
\begin{array}{ccc}
a_{LN}\left(k\right) & = & \sqrt{\pi}(2L+1)2^{-L-1}i^{L-\text{N}}\sum_{M=0}^{\infty}\frac{\left(\left(-\frac{1}{4}\right)^{M}k^{L+2M}\right)}{2^{L+2M+1}\left(M!\Gamma\left(L+M+\frac{3}{2}\right)\right)}\\
 & \times & \left(1+(-1)^{L+2M+\text{N}}\right)\binom{L+2M}{\frac{1}{2}(L+2M-\text{N})}\\
 & = & \frac{\sqrt{\pi}2^{-2L-2}(2L+1)k^{L}i^{L-N}}{\Gamma\left(\frac{1}{2}(2L+3)\right)}\left(1+(-1)^{L+N}\right)\binom{L}{\frac{L-N}{2}}\\
 & \times & \,_{2}F_{3}\left(\frac{L}{2}+\frac{1}{2},\frac{L}{2}+1;L+\frac{3}{2},\frac{L}{2}-\frac{N}{2}+1,\frac{L}{2}+\frac{N}{2}+1;-\frac{k^{2}}{4}\right)\\
 & = & \sqrt{\pi}\,2^{-2L-2}(2L+1)k^{L}i^{L-N}\left(1+(-1)^{L+N}\right)\Gamma(L+1)\\
 & \times & \,_{2}\tilde{F}_{3}\left(\frac{L}{2}+\frac{1}{2},\frac{L}{2}+1;L+\frac{3}{2},\frac{L}{2}-\frac{N}{2}+1,\frac{L}{2}+\frac{N}{2}+1;-\frac{k^{2}}{4}\right)\end{array}\;,\label{eq:a_L_as_2F3}\end{equation}
of which the final two steps were new with the prior work \cite{stra24b}.
We included the final form using regularized hypergeometric functions
\cite{wolfram.com/07.26.26.0001.01}

\begin{equation}
\,_{2}F_{3}\left(a_{1},a_{2};b_{1},b_{2},b_{3};z\right)=\Gamma\left(b_{1}\right)\Gamma\left(b_{2}\right)\Gamma\left(b_{3}\right)\,_{2}\tilde{F}_{3}\left(a_{1},a_{2};b_{1},b_{2},b_{3};z\right)\label{eq:regularized}\end{equation}
and cancelled the $\Gamma\left(b_{i}\right)$ with gamma functions
in the denominators of the prefactors. This cancellation allows one
to avoid infinities that arise whenever \emph{N>1} is an integer larger
than \emph{L}, and of the same parity, that would otherwise result in indeterminacies
in a computation when one tries to use the conventional form of the
hypergeometric function.

After a further review of the literature we found that Keatings result
(the first line above) and ours (the second line above) are implicitly
subsumed within Jet Wimp's 1962 Jacobi expansion \cite{jacobi_cheb_Wimp}of
the Anger-Weber function (his equations (2.10) and (2.11)) since Legendre
polynomials are a subset of Jacobi polynomials and the Bessel function
$J_{N}\left(kx\right)$ are a special case of the Anger-Weber function
$\mathbf{J}_{\nu}\left(kx\right)$ when $\nu$ is an integer. Wimp
does not mention the calculational difficulties that were resolved
through the third form, above.

For the special cases of $N=0,\,1$ the order of the hypergeometric
functions is reduced since the parameters $a_{2}=b_{3}$ and $a_{1}=b_{2}$,
resp., giving

\begin{equation}
\begin{array}{ccc}
a_{L0}\left(k\right) & = & \frac{\sqrt{\pi}i^{L}2^{-2L-2}(2L+1)k^{L}}{\Gamma\left(\frac{1}{2}(2L+3)\right)}\left(1+(-1)^{L}\right)\binom{L}{\frac{L}{2}}\\
 & \times & \,_{1}F_{2}\left(\frac{L}{2}+\frac{1}{2};\frac{L}{2}+1,L+\frac{3}{2};-\frac{k^{2}}{4}\right)\\
 & = & \sqrt{\pi}i^{L}2^{-2L-2}(2L+1)k^{L}\Gamma\left(\frac{L}{2}+1\right)\left(1+(-1)^{L}\right)\binom{L}{\frac{L}{2}}\\
 & \times & \,_{1}\tilde{F}_{2}\left(\frac{L}{2}+\frac{1}{2};\frac{L}{2}+1,L+\frac{3}{2};-\frac{k^{2}}{4}\right)\end{array}\;,\label{eq:a_L0}\end{equation}

and

\begin{equation}
\begin{array}{ccc}
a_{L1}\left(k\right) & = & \frac{\sqrt{\pi}i^{L-1}2^{-2L-2}(2L+1)k^{L}}{\Gamma\left(\frac{1}{2}(2L+3)\right)}\left(1+(-1)^{L+1}\right)\binom{L}{\frac{L-1}{2}}\\
 & \times & \,_{1}F_{2}\left(\frac{L}{2}+1;\frac{L}{2}+\frac{3}{2},L+\frac{3}{2};-\frac{k^{2}}{4}\right)\\
 & = & i^{L-1}2^{-L-2}(2L+1)k^{L}\Gamma\left(\frac{L}{2}+1\right)\left(1+(-1)^{L+1}\right)\\
 & \times & \,_{1}\tilde{F}_{2}\left(\frac{L}{2}+1;\frac{L}{2}+\frac{3}{2},L+\frac{3}{2};-\frac{k^{2}}{4}\right)\end{array}\;,\label{eq:a_L1}\end{equation}
In each special case the first form involving a hypergeometric function
has no indeterminacies, but we include the regularized hypergeometric
function version for completeness.

The first 22 terms in the Fourier-Legendre series for $J_{0}\left(kx\right)$
(\ref{eq:Fourier-Legendre series}) are then, with $k=1$ \begin{eqnarray}
J_{0}\left(x\right) & \cong & 0.9197304100897602393144211940806200P_{0}(x)\nonumber \\
 & - & 0.1579420586258518875737139671443637P_{2}(x)\nonumber \\
 & + & 0.003438400944601109232996887872072915P_{4}(x)\nonumber \\
 & - & 0.00002919721848828729693660590986125663P_{6}(x)\nonumber \\
 & + & 1.317356952447780977655616563143280\times10^{-7}\; P_{8}(x)\nonumber \\
 & - & 3.684500844208203027173771096058866\times10^{-10}\; P_{10}(x)\nonumber \\
 & + & 7.011830032993845928208803328211457\times10^{-13}\; P_{12}(x)\nonumber \\
 & - & 9.665964369858912263671995372753346\times10^{-16}\; P_{14}(x)\nonumber \\
 & + & 1.009636276824546446525342170924936\times10^{-18}\; P_{16}(x)\nonumber \\
 & - & 8.266656955927637858991972584174117\times10^{-22}\; P_{18}(x)\nonumber \\
 & + & 5.448244867762758725890082837839430\times10^{-25}\; P_{20}(x)\nonumber \\
 & - & 2.952527182137354751675774606663400\times10^{-28}\; P_{22}(x)\nonumber \\
 & + & 1.338856158858534469080898670096200\times10^{-31}\; P_{24}(x)\nonumber \\
 & - & 5.154913186088512926193234837816582\times10^{-35}\; P_{26}(x)\nonumber \\
 & + & 1.706231577038503450138564028467634\times10^{-38}\; P_{28}(x)\nonumber \\
 & - & 4.906893556427796857473097979568289\times10^{-42}\; P_{30}(x)\nonumber \\
 & + & 1.237489200717479383020539576221293\times10^{-45}\; P_{32}(x)\nonumber \\
 & - & 2.759056237537871868604555688548364\times10^{-49}\; P_{34}(x)\nonumber \\
 & + & 5.477382207172712629199714648396409\times10^{-53}\; P_{36}(x)\nonumber \\
 & - & 9.744200345578852550688946057050674\times10^{-57}\; P_{38}(x)\nonumber \\
 & + & 1.562280711659504489828025148995770\times10^{-60}\; P_{40}(x)\nonumber \\
 & - & 2.269056283827394368836057470594599\times10^{-64}\; P_{42}(x)\;.\label{eq:Fourier-Legendre series-numerical}\end{eqnarray}

\noindent
In equation (15) of  the prior work \cite{stra24b} the last line, above, mistakenly had the wrong power, 

\noindent
$-2.269056283827394368836057470594599\times10^{-60}\; P_{42}(x)$, though the Fortran code in the appendix was correct: -2.269056283827394368836057470594599 e-64 P(42 , x ). 

One can easily expand the Legendre polynomials into their constituent
terms and gather like powers in (\ref{eq:Fourier-Legendre series-numerical})
to give an updated polynomial approximation, \begin{eqnarray}
J_{0}\left(x\right) & \cong & 1.000000000000000000000000000000000000\, x^{0}\nonumber \\
 & - & 0.2500000000000000000000000000000000000\, x^{2}\nonumber \\
 & + & 0.01562500000000000000000000000000000000\, x^{4}\nonumber \\
 & - & 0.0004340277777777777777777777777777777778\, x^{6}\nonumber \\
 & + & 6.781684027777777777777777777777777778\times10^{-6}\; x^{8}\nonumber \\
 & - & 6.781684027777777777777777777777777778\times10^{-8}\; x^{10}\nonumber \\
 & + & 4.709502797067901234567901234567901235\times10^{-10}\; x^{12}\nonumber \\
 & - & 2.402807549524439405391786344167296548\times10^{-12}\; x^{14}\nonumber \\
 & + & 9.385966990329841427311665406903502142\times10^{-15}\; x^{16}\nonumber \\
 & - & 2.896903392077111551639402903365278439\times10^{-17}\; x^{18}\nonumber \\
 & + & 7.242258480192778879098507258413196097\times10^{-20}\; x^{20}\nonumber \\
 & - & 1.496334396734045222954237036862230599\times10^{-22}\; x^{22}\nonumber \\
 & + & 2.597802772107717400962217077885817011\times10^{-25}\; x^{24}\nonumber \\
 & - & 3.842903509035084912666001594505646466\times10^{-28}\; x^{26}\nonumber \\
 & + & 4.901662639075363409012757135849038860\times10^{-31}\; x^{28}\nonumber \\
 & - & 5.446291821194848232236396817610043178\times10^{-34}\; x^{30}\nonumber \\
 & + & 5.318644356635593976793356267197307791\times10^{-37}\; x^{32}\nonumber \\
 & - & 4.600903422695150498956190542558224733\times10^{-40}\; x^{34}\nonumber \\
 & + & 3.550079801462307483762492702591222788\times10^{-43}\; x^{36}\nonumber \\
 & - & 2.458504017633176927813360597362342651\times10^{-46}\; x^{38}\nonumber \\
 & + & 1.5365650110207355798833503733514641567\times10^{-49}\; x^{40}\nonumber \\
 & - & 8.7106860035189091830121903251216788929\times10^{-53}\; x^{42}\nonumber \\
 & \cong & 1-\frac{x^{2}}{2^{2}}+\frac{x^{4}}{2^{6}}-\frac{x^{6}}{2^{8}3^{2}}+\frac{x^{8}}{2^{14}3^{2}}-\frac{x^{10}}{2^{16}3^{2}5^{2}}+\frac{x^{12}}{2^{20}3^{4}5^{2}}-\frac{x^{14}}{2^{22}3^{4}5^{2}7^{2}}+\frac{x^{16}}{2^{30}3^{4}5^{2}7^{2}}\nonumber \\
 & - & \frac{x^{18}}{2^{32}3^{8}5^{2}7^{2}}+\frac{x^{20}}{2^{36}3^{8}5^{4}7^{2}}-\frac{x^{22}}{2^{38}3^{8}5^{4}7^{2}11^{2}}+\frac{x^{24}}{2^{44}3^{10}5^{4}7^{2}11^{2}}-\frac{x^{26}}{2^{46}3^{10}5^{4}7^{2}11^{2}13^{2}}\nonumber \\
 & + & \frac{x^{28}}{2^{50}3^{10}5^{4}7^{4}11^{2}13^{2}}-\frac{x^{30}}{2^{52}3^{12}5^{6}7^{4}11^{2}13^{2}}+\frac{x^{32}}{2^{62}3^{12}5^{6}7^{4}11^{2}13^{2}}\nonumber \\
 & - & \frac{x^{34}}{2^{64}3^{12}5^{6}7^{4}11^{2}13^{2}17^{2}}+\frac{x^{36}}{2^{68}3^{16}5^{6}7^{4}11^{2}13^{2}17^{2}}-\frac{x^{38}}{2^{70}3^{16}5^{6}7^{4}11^{2}13^{2}17^{2}19^{2}}\nonumber \\
 & + & \frac{x^{40}}{2^{76}3^{16}5^{8}7^{4}11^{2}13^{2}17^{2}19^{2}}-\frac{x^{42}}{2^{78}3^{18}5^{8}7^{6}11^{2}13^{2}17^{2}19^{2}}\;.\label{eq:poly_approx_48_digit}\end{eqnarray}
 The latter form is simply an inverse prime version of the first 22
terms of the well-known series representation~\cite{GR5 p. 987 No. 8.511.4}
(p. 970 No. 8.440) \begin{equation}
J_{\nu}\left(x\right)=\sum_{k=0}^{\infty}\frac{(-1)^{k}\left(\frac{x}{2}\right)^{2k+\nu}}{k!\Gamma(k+\nu+1)}.\label{eq:GR5 p. 970 No. 8.440}\end{equation}
Since each Legendre polynomial contributes to the constant term in
both (\ref{eq:poly_approx_48_digit}) and (\ref{eq:GR5 p. 970 No. 8.440})
their sum is 

\begin{eqnarray}
0.919730410089760239314421194080620\times & (1)\nonumber \\
-0.157942058625851887573713967144364\times & \frac{1}{2}\left(3x^{2}-1\right)_{x\rightarrow0}\nonumber \\
+0.00343840094460110923299688787207292\times & \frac{1}{8}\left(35x^{4}-30x^{2}+3\right)_{x\rightarrow0}\nonumber \\
-0.0000291972184882872969366059098612566\times & \frac{1}{16}\left(231x^{6}-315x^{4}+105x^{2}-5\right)_{x\rightarrow0}\nonumber \\
+\cdots & = & 1\label{eq:first4products}\end{eqnarray}
rather than some other number close to \emph{1}. This may formalized
in the summed series 

\begin{eqnarray}
\sum_{L=0}^{\infty}\, &  & \hspace{-0.9cm}\frac{\sqrt{\pi}i^{L}2^{-2L-2}\left(1+(-1)^{L}\right)(2L+1)\binom{L}{\frac{L}{2}}}{\Gamma\left(\frac{1}{2}(2L+3)\right)}\,_{1}F_{2}\left(\frac{L}{2}+\frac{1}{2};\frac{L}{2}+1,L+\frac{3}{2};-\frac{k^{2}}{4}\right)k^{L}\nonumber \\
 & \times & \left\{ \frac{i^{L}2^{-L/2}(L-1)\text{!!}\left(\frac{L}{2}+\frac{1}{2}\right)_{h}\left(-\frac{L}{2}\right)_{h}}{h!\frac{L}{2}!\left(\frac{1}{2}\right)_{h}}\right\} =\frac{(-1)^{h}2^{-2h}}{h!\Gamma(h+1)}k^{2h}\;,\label{eq:1/primes^n=00003D00003D00003D}\end{eqnarray}
within which $h=0$ gives (\ref{eq:first4products}). This may be
derived for the general \emph{h} by extracting specific powers of
\emph{x} from the Legendre polynomials, most easily by converting
them into $\,_{2}F_{1}$ hypergeometric functions \cite{GR5 p. 987 No. 8.511.4} (p. 1044 No.
8911.2)  \cite{PBM3 p. 468 No. 7.3.1.206}
(p. 466 No. 7.3.1.182) and thence into a finite sum over ratios of Pochhammer symbols.

In equation (26) of the previous paper \cite{stra24b} the term in curly braces was
given as

\begin{equation}
\left\{ \frac{2^{-L}\binom{2L}{L}\left(\frac{1}{2}-\frac{L}{2}\right)_{\frac{L}{2}-h}\left(-\frac{L}{2}\right)_{\frac{L}{2}-h}}{\left(\frac{L}{2}-h\right)!\left(\frac{1}{2}-L\right)_{\frac{L}{2}-h}}\right\} \label{eq:orig_J0_series}\end{equation}
because it used an alternaticve conversion of Legendre polynomials
into $\,_{2}F_{1}$ hypergeometric functions \cite{GR5 p. 987 No. 8.511.4} (p. 1044 No.
8911.1)  \cite{PBM3 p. 468 No. 7.3.1.206}
(p. 468 No. 7.3.1.206).

For $h=1$, the $P_{2}(x)$ through $P_{42}(x)$ terms all add to
give $-1/4$, the coefficient of $x^{2}$ term in both (\ref{eq:poly_approx_48_digit})
and (\ref{eq:GR5 p. 970 No. 8.440}) if $k=1$. Including $k\neq1$
poses no problem in (\ref{eq:1/primes^n=00003D00003D00003D}) despite
its appearance as the argumemt of the $\,_{1}F_{2}\left(\frac{L}{2}+\frac{1}{2};\frac{L}{2}+1,L+\frac{3}{2};-\frac{k^{2}}{4}\right)$
function, as well as there being a $k^{L}$ factor in the argument
of the sum. It ends up contributing a very clean factor of $k^{2h}$
to the right-hand side of (\ref{eq:1/primes^n=00003D00003D00003D}).
We thus have summed an infinite set of infinte sums of $\,_{1}F_{2}$
hypergeometric functions, though we have verified only those with
$0\leq h\leq42$. (We had to take the upper limit on the number of
terms in the series $\geq h+74$ in order to obtain a percent difference
between left- and right-hand sides that was $\leq10^{-33}$ because
the first \emph{h} terms in the series do not contribute. For $h=0$,
an upper limit on the number of terms in the series $\geq h+44$ was
sufficient.)

We likewise summed an infinite set of infinte sums of $\,_{1}F_{2}$
hypergeometric functions derived from the Fourier-Legendre series
for $J_{1}\left(kx\right)$ (\ref{eq:Fourier-Legendre series}): \begin{eqnarray}
\sum_{L=1}^{\infty}\, &  & \hspace{-0.9cm}\frac{\sqrt{\pi}i^{L-1}\left(1+(-1)^{L+1}\right)(2L+1)2^{-3L-2}\binom{L}{\frac{L-1}{2}}\binom{2L}{L}(-1)^{-h+\frac{L}{2}-\frac{1}{2}}}{\Gamma\left(\frac{1}{2}(2L+3)\right)\left(-h+\frac{L}{2}-\frac{1}{2}\right)!}k^{L}\nonumber \\
 & \times & \,_{1}F_{2}\left(\frac{L}{2}+1;\frac{L}{2}+\frac{3}{2},L+\frac{3}{2};-\frac{k^{2}}{4}\right)\left[\frac{\left(\frac{1}{2}-\frac{L}{2}\right)_{\frac{L-1}{2}-h}\left(-\frac{L}{2}\right)_{\frac{L-1}{2}-h}}{\left(\frac{1}{2}-L\right)_{\frac{L-1}{2}-h}}\right]\nonumber \\
 & = & \frac{(-1)^{h}2^{-2h-1}}{h!\Gamma(h+2)}k^{2h+1}\;.\label{eq:J1summedtoprimes}\end{eqnarray}

The question natually arises as to whether one can derived such summed
infinite series based on other polynomial expansions. In the following,
one may answer in the affirmative for both Chebyshev and Gegenbauer
polynomial expansions of Bessel functions.

\section{Summed series involving $\,_{1}F_{2}$ hypergeometric functions from
Chebyshev polynomial expansions of Bessel functions}

Wimp also applied his Jacobi expansion \cite{jacobi_cheb_Wimp}to
find Chebyshev polynomial expansions of Bessel functions, since \cite{GR5 p. 987 No. 8.511.4}
(p. 1060 No. 8.962.3)\begin{equation}
P_{2n}^{\left(-\frac{1}{2},-\frac{1}{2}\right)}(z)=\frac{\left(\frac{1}{2}\right)_{2n}}{(2n)!}T_{2n}(z)\;.\label{eq:Jacobi_Chebyshev}\end{equation}
Unlike the section above, the following expansion (his Equations (3.6)
and (3.7)) applies to non-integer indices as well:

\begin{equation}
J_{\upsilon}\left(kx\right)=\left(kx\right)^{\upsilon}\sum_{L=0}^{\infty}C_{L\upsilon}\left(k\right)T_{2L}(x)\:.\quad-1\leq x\leq1\label{eq:Chebyshev series}\end{equation}
Since what we are expanding in Chebyshev polynomials is the function
$J_{\upsilon}\left(kx\right)\left(kx\right)^{-\upsilon}$, the coefficients
can only be given by the orthogonality of the Chebyshev polynomials
if we include the full function in the integral, \begin{equation}
C_{L\upsilon}\left(k\right)=\frac{\left(2-\delta_{L0}\right)}{2}\frac{2}{\pi}\int_{-1}^{1}\left(J_{\upsilon}\left(kx\right)\left(kx\right)^{-\upsilon}\right)\frac{1}{\sqrt{1-x^{2}}}T_{2L}(x)\, dx\;,\label{eq:coefficients-1}\end{equation}
which Wimp finds to be

\begin{equation}
\begin{array}{ccc}
C_{L\upsilon}\left(k\right) & = & {\displaystyle \frac{(-1)^{L}k^{2L}2^{-4L-\nu}\left(2-\delta_{L0}\right)}{L!\Gamma(L+\nu+1)}}\,_{1}F_{2}\left(L+\frac{1}{2};L+\nu+1,2L+1;-\frac{k^{2}}{4}\right)\end{array}\;,\label{eq:C_L_as_1F2}\end{equation}

The first 22 terms in the Chebyshev polynomial expansions of $J_{0}\left(kx\right)$
(\ref{eq:Fourier-Legendre series}) are then, with $k=1$ \begin{eqnarray}
J_{0}\left(x\right) & \cong & 0.8807255791026085285666716907449594\, T_{0}(x)\nonumber \\
 & - & 0.1173880111683243194062454639255572\, T_{2}(x)\nonumber \\
 & + & 0.001873212503719194837870878203929524\, T_{4}(x)\nonumber \\
 & - & 0.00001314542297029262107182993119503582\, T_{6}(x)\nonumber \\
 & + & 5.167242966801437053171032359951600\times10^{-8}\; T_{8}(x)\nonumber \\
 & - & 1.297218234854703963093975334759865\times10^{-10}\; T_{10}(x)\nonumber \\
 & + & 2.258840234607001930320227243984034\times10^{-13}\; T_{12}(x)\nonumber \\
 & - & 2.887621352768057764464058481597816\times10^{-16}\; T_{14}(x)\nonumber \\
 & + & 2.824848256251380023621233536051211\times10^{-19}\; T_{16}(x)\nonumber \\
 & - & 2.182699061309088513825726048290021\times10^{-22}\; T_{18}(x)\nonumber \\
 & + & 1.365739183823366078819378297317202\times10^{-25}\; T_{20}(x)\nonumber \\
 & - & 7.061125701699520180896051661348297\times10^{-29}\; T_{22}(x)\nonumber \\
 & + & 3.067182727248138051740188483703613\times10^{-32}\; T_{24}(x)\nonumber \\
 & - & 1.135092833714987500414966932525964\times10^{-35}\; T_{26}(x)\nonumber \\
 & + & 3.621712251769489873248477093327996\times10^{-39}\; T_{28}(x)\nonumber \\
 & - & 1.006555480914216913705134524512148\times10^{-42}\; T_{30}(x)\nonumber \\
 & + & 2.458540787185135207907001122952213\times10^{-46}\; T_{32}(x)\nonumber \\
 & - & 5.319086471776732419423425079488687\times10^{-50}\; T_{34}(x)\nonumber \\
 & + & 1.026433533066142649943339190916424\times10^{-53}\; T_{36}(x)\nonumber \\
 & - & 1.777651158721406916387585852076982\times10^{-57}\; T_{38}(x)\nonumber \\
 & + & 2.778406892667094352173643013096289\times10^{-61}\; T_{40}(x)\nonumber \\
 & - & 3.938717221679009654181092747102998\times10^{-65}\; T_{42}(x)\quad-1\leq x\leq1\;.\label{eq:J0 in Chebyshev}\end{eqnarray}

At the upper limit of applicability, $x=1$, this gives 33-digit accuracy,
$J_{0}\left(1\right)=0.765197686557966551449717526102663$. (Even
at $x=8$, this gives a result accurate to 14-digits, $J_{0}\left(8\right)=0.171650807137554$.) 

If one follows Clenshaw's \cite{Clenshaw} (p. 30) lead and instead
takes $k=8$ one gets

\begin{eqnarray}
J_{0}\left(x\right) & \cong & \mathbf{0.315455942949780239127}5502330199159/2\, T_{0}(x/8)\nonumber \\
 & - & \mathbf{0.00872344235285222129}0793322469895429\, T_{2}(x/8)\nonumber \\
 & + & \mathbf{0.265178613203336809867}0778235911043\, T_{4}(x/8)\nonumber \\
 & - & \mathbf{0.37009499387264977903}34193036836753\, T_{6}(x/8)\nonumber \\
 & + & \mathbf{0.15806710233209726127}77155496720475\; T_{8}(x/8)\nonumber \\
 & - & \mathbf{0.03489376941140888516}317328987958171\; T_{10}(x/8)\nonumber \\
 & + & \mathbf{0.00481918006946760449}6778380314312767\; T_{12}(x/8)\nonumber \\
 & - & \boldsymbol{0.000460626166206275047}5036418408154116\; T_{14}(x/8)\nonumber \\
 & + & \mathbf{0.00003246032882100508}080625560924485746\; T_{16}(x/8)\nonumber \\
 & - & \mathbf{1.76194690776215}0749459765966407618\times10^{-6}\; T_{18}(x/8)\nonumber \\
 & + & \mathbf{7.6081635924187}81866973786230699492\times10^{-8}\; T_{20}(x/8)\nonumber \\
 & - & \mathbf{2.679253530557}672898335371633826306\times10^{-9}\; T_{22}(x/8)\nonumber \\
 & + & \mathbf{7.848696314}479464416529503905101749\times10^{-11}\; T_{24}(x/8)\nonumber \\
 & - & \mathbf{1.943834686}737016570620688424557753\times10^{-12}\; T_{26}(x/8)\nonumber \\
 & + & \mathbf{4.1253205}95634373932612618412757652\times10^{-14}\; T_{28}(x/8)\nonumber \\
 & - & \mathbf{7.5885}08125447546337619860819329317\times10^{-16}\; T_{30}(x/8)\nonumber \\
 & + & \mathbf{1.2218}51587396141103441861977201729\times10^{-17}\; T_{32}(x/8)\nonumber \\
 & - & \mathbf{1.7}36789607700236768293730242713393\times10^{-19}\; T_{34}(x/8)\nonumber \\
 & + & 2.195793203319560353679493897698779\times10^{-21}\; T_{36}(x/8)\nonumber \\
 & - & 2.485566419364292266537947175258836\times10^{-23}\; T_{38}(x/8)\nonumber \\
 & + & 2.534024606818972691102585769070259\times10^{-25}\; T_{40}(x/8)\nonumber \\
 & - & 2.339085627055744706712023052059754\times10^{-28}\; T_{42}(x/8)\quad-8\leq x\leq8\;.\label{eq:J0 in Chebyshev x/8}\end{eqnarray}
where the bolding indicates the digits he displays. (I have included
an extra digit in some places to allow for appropriate rounding to
his displayed digit.) Clenshaw follows the usual convention (noted
on his p. 7) for sums having a single prime to indicate that the term
with suffix zero is to be halved (and if the prime is doubled, the
highest term in the sum is also halved), as indicated in the first
line of (\ref{eq:J0 in Chebyshev x/8}). This factor of two difference
arises from the normalization of the orthogonality relation for Chebyshev
polynomials \cite{GR5 p. 987 No. 8.511.4} (p. 1057 No. 8.949.9):

\begin{equation}
\int_{-1}^{1}T_{n}(x)T_{m}(x)\frac{1}{\sqrt{1-x^{2}}}=\begin{cases}
\begin{array}{c}
0\\
\pi/2\\
\pi\end{array} & \begin{array}{c}
m\neq n\\
m=n\neq0\\
m=n=0\end{array}\end{cases}\:.\label{eq:T_orth}\end{equation}

Since I am comparing Chebyshev expansions with both Legendre and Gegenbauer
expansions, whose orthogonality relations \cite{GR5 p. 987 No. 8.511.4}
(p. 1043 No. 8.904 and p. 1054 No. 8.939.8, resp.) have no such third
branch, all derivations are made much more straightforward if one
adopts the perhaps iconclastic notion that sums should simply be sums
and simply displays the first line of (\ref{eq:J0 in Chebyshev x/8})
as $0.1577279714748901195637751165099580\, T_{0}(x/8)$. 

At the upper limit of applicability, $x=8$, (\ref{eq:J0 in Chebyshev x/8})
gives 27-digit accuracy, $J_{0}\left(8\right)=0.171650807137553906090869408$. 

Since the constant term of every Chebyshev polynomial has magnitude
one, and alternating sign, the sum of these times the coefficients
--$a_{2r}$ in Clenshaw's convention -- is simply 

\begin{equation}
\sum\,'\left(-1\right)^{L}a_{2r}=1\:.\label{eq:Clenshaw sum}\end{equation}
Clenshaw has, thus, given the first of the summation rules we wish
to derive.

In order to sum the infinite set of infinte sums of $\,_{1}F_{2}$
hypergeometric functions derived from the Fourier-Legendre series
expansion of Bessel functions (\ref{eq:1/primes^n=00003D00003D00003D})
and (\ref{eq:J1summedtoprimes}), we extracted specific powers of
\emph{x} by converting Legendre polynomials into $\,_{2}F_{1}$ hypergeometric
functions \cite{PBM3 p. 468 No. 7.3.1.206} (p. 468 No. 7.3.1.206)
and thence into a finite sum over ratios of Pochhammer symbols. The
equivalent conversion for Chebyshev polynomials is

\begin{equation}
T_{2L}(x)=2^{2L-1}x^{2L}\,_{2}F_{1}\left(\frac{1}{2}-L,-L;1-2L;\frac{1}{x^{2}}\right)\left[1+\delta_{L0}\right]\:.\label{eq:PBM3 p. 468 No. 7.3.1.207}\end{equation}

Note that I have augmented \cite{PBM3 p. 468 No. 7.3.1.206} (p. 468
No. 7.3.1.207) \cite{functions.wolfram.com/07.23.03.0191.01} with
the factor $\left[1+\delta_{L0}\right]$ that allows the conversion
to be extanded downward from $L>0$. When multiplied by the equivalent
factor in Wimp's Chebyshev expansion (\ref{eq:C_L_as_1F2}), one obtains
$\left(2-\delta_{L0}\right)\left[1+\delta_{L0}\right]\equiv2$ for
all \emph{L}. This is a strong argument for using the {}``sums should
simply be sums'' convention over Clenshaw's for the present analytical
work, so that the summed series associated with the power $x^{2h}$
may simply be written as

\begin{eqnarray}
\sum_{L=0}^{\infty}\, &  & \hspace{-0.9cm}\frac{(-1)^{L}2^{-2L}\left(\frac{1}{2}-L\right)_{L-h}(-L)_{L-h}}{L!\Gamma(L+1)(L-h)!(1-2L)_{L-h}}\,_{1}F_{2}\left(L+\frac{1}{2};L+1,2L+1;-\frac{k^{2}}{4}\right)k^{2L}\nonumber \\
 & = & \frac{(-1)^{h}2^{-2h}}{h!\Gamma(h+1)}k^{2h}\;.\label{eq:J0_T_series}\end{eqnarray}
To verify the lowest $43$ summed series for $k\rightarrow8$, one
has to take the upper limit on the number of terms in the series $\geq h+18$
in order to obtain a percent difference between left- and right-hand
sides that was $\leq10^{-33}$ because the first \emph{h} terms in
the series do not contribute. For $k\rightarrow5$, this reduces somewhat
to $\geq h+15$. For $h=0$, one needs 24 terms and 20 terms, respectively. 

The first 22 terms in the $J_{1}\left(x\right)$ Chebyshev expansion
(\ref{eq:Chebyshev series}) with $k=1$ are \begin{eqnarray}
J_{1}\left(x\right) & \cong & 0.4697097923433853441348972113538690xT_{0}(x)\nonumber \\
 & - & 0.02997305358809894507094444118401190xT_{2}(x)\nonumber \\
 & + & 0.0003154953401761330198307113032804328xT_{4}(x)\nonumber \\
 & - & 1.653528591827665010389921139509211\times10^{-6}\; xT_{6}(x)\nonumber \\
 & + & 5.188889110114106792954599573058750\times10^{-9}\; xT_{8}(x)\nonumber \\
 & - & 1.084245120515337519078432469943857\times10^{-11}\; xT_{10}(x)\nonumber \\
 & + & 1.617069529094057869823401928778476\times10^{-14}\; xT_{12}(x)\nonumber \\
 & - & 1.807903976592524723392831520195131\times10^{-17}\; xT_{14}(x)\nonumber \\
 & + & 1.571543945521723529179083698815771\times10^{-20}\; xT_{16}(x)\nonumber \\
 & - & 1.092591641508275242057122355553840\times10^{-23}\; xT_{18}(x)\nonumber \\
 & + & 6.213791797992245609440469557453575\times10^{-27}\; xT_{20}(x)\nonumber \\
 & - & 2.944495823790016197177000782247634\times10^{-30}\; xT_{22}(x)\nonumber \\
 & + & 1.180496667850251944095467073781979\times10^{-33}\; xT_{24}(x)\nonumber \\
 & - & 4.056318036675064198378921654189439\times10^{-37}\; xT_{26}(x)\nonumber \\
 & + & 1.207866649436639014639549760562102\times10^{-40}\; xT_{28}(x)\nonumber \\
 & - & 3.146932355403406273096620834992699\times10^{-44}\; xT_{30}(x)\nonumber \\
 & + & 7.233957871819338833114752440681911\times10^{-48}\; xP_{32}(x)\nonumber \\
 & - & 1.478064332069756593976138661523809\times10^{-51}\; xT_{34}(x)\nonumber \\
 & + & 2.702029827426988943325772959142285\times10^{-55}\; xT_{36}(x)\nonumber \\
 & - & 4.445451117805773022660415901032200\times10^{-59}\; xT_{38}(x)\nonumber \\
 & + & 6.617045043041664246398527226007578\times10^{-63}\; xT_{40}(x)\nonumber \\
 & - & 8.953842205918258708007813804592169\times10^{-67}\; xT_{42}(x)\quad-1\leq x\leq1\;.\label{eq:Chebyshev_J1}\end{eqnarray}
At the upper limit of applicability, $x=1$, this gives 33-digit accuracy,
$J_{1}\left(1\right)=0.440050585744933515959682203718915$. (Even
at $x=8$, this gives a result accurate to 16-digits, $J_{1}\left(8\right)=0.2346363468539146$.)

If one follows Clenshaw's \cite{Clenshaw} (p. 31) lead and instead
takes $k=8$ one gets

\begin{eqnarray}
J_{1}\left(x\right) & \cong & \mathbf{1.29671754121052984167}3374221959245/2\frac{x}{8}T_{0}(\frac{x}{8})\nonumber \\
 & - & \mathbf{1.191801160541216872507}032741866674\frac{x}{8}T_{2}(\frac{x}{8})\nonumber \\
 & + & \mathbf{1.28799409885767762038}2580899489350\frac{x}{8}T_{4}(\frac{x}{8})\nonumber \\
 & - & \mathbf{0.66144393413454325277}28770946844658\frac{x}{8}T_{6}(\frac{x}{8})\nonumber \\
 & + & \mathbf{0.17770911723972828328}23229884383241\frac{x}{8}T(\frac{x}{8})\nonumber \\
 & - & \mathbf{0.02917552480615420766}201489599627591\frac{x}{8}T_{10}(\frac{x}{8})\nonumber \\
 & + & \mathbf{0.00324027018268385746}6456539040415511\frac{x}{8}T_{12}(\frac{x}{8})\nonumber \\
 & - & \mathbf{0.00026044438934858068}13446141103993105\frac{x}{8}T_{14}(\frac{x}{8})\nonumber \\
 & + & \mathbf{0.00001588701923993213}39310461547076296\frac{x}{8}T_{16}(\frac{x}{8})\nonumber \\
 & - & \mathbf{7.6175878054003}48945692364404508548\times10^{-7}\;\frac{x}{8}T_{18}(\frac{x}{8})\nonumber \\
 & + & \mathbf{2.9497070072777}18590826100996112190\times10^{-8}\;\frac{x}{8}T_{20}(\frac{x}{8})\nonumber \\
 & - & \mathbf{9.42421298156}7078718578173809056009\times10^{-10}\;\frac{x}{8}T_{22}(\frac{x}{8})\nonumber \\
 & + & \mathbf{2.528123664}278402657192198903253796\times10^{-11}\;\frac{x}{8}T_{24}(\frac{x}{8})\nonumber \\
 & - & \mathbf{5.77740419}1721418742769122933910453\times10^{-13}\;\frac{x}{8}T_{26}(\frac{x}{8})\nonumber \\
 & + & \mathbf{1.1385715}20281115385303951328717824\times10^{-14}\;\frac{x}{8}T_{28}(\frac{x}{8})\nonumber \\
 & - & \mathbf{1.95535}7833295237111457156049739834\times10^{-16}\;\frac{x}{8}T_{30}(\frac{x}{8})\nonumber \\
 & + & \mathbf{2.95}3014639834346609722058184262545\times10^{-18}\;\frac{x}{8}P_{32}(\frac{x}{8})\nonumber \\
 & - & \mathbf{3.9}52934614113459501768862170679755\times10^{-20}\;\frac{x}{8}T_{34}(\frac{x}{8})\nonumber \\
 & + & 4.723067439441036227167716497766825\times10^{-22}\;\frac{x}{8}T_{36}(\frac{x}{8})\nonumber \\
 & - & 5.068481382508651457731548219527637\times10^{-24}\;\frac{x}{8}T_{38}(\frac{x}{8})\nonumber \\
 & + & 4.912426488809207456168647750374833\times10^{-26}\;\frac{x}{8}T_{40}(\frac{x}{8})\nonumber \\
 & - & 4.321688707060755263766813871186111\times10^{-28}\;\frac{x}{8}T_{42}(\frac{x}{8})\quad-8\leq x\leq8\;.\label{eq:Clenshaw_J1_k=00003D8}\end{eqnarray}
where the bolding indicates the digits he displays. (I have included
an extra digit in some places to allow for appropriate rounding to
his displayed digit.) If one takes the iconoclastic route of not followed
his convention (noted on his p. 7) for sums having a single prime
to indicate that the term with suffix zero is to be halved, the first
line, above would be $0.6483587706052649208366871109796227\,\frac{x}{8}T_{0}(\frac{x}{8})$. 

At the upper limit of applicability, $x=8$, (\ref{eq:Clenshaw_J1_k=00003D8})
gives 29-digit accuracy, $J_{1}\left(8\right)=0.23463634685391462438127665159$. 

The summed series associated with associated with the power $x^{2h+1}$
in this Chebyshev expansion are 

\begin{eqnarray}
\sum_{L=1}^{\infty}\, &  & \hspace{-0.9cm}\frac{(-1)^{L}2^{-2L-1}\left(\frac{1}{2}-L\right)_{L-h}(-L)_{L-h}}{L!\Gamma(L+2)(L-h)!(1-2L)_{L-h}}k^{2L+1}\,_{1}F_{2}\left(L+\frac{1}{2};L+2,2L+1;-\frac{k^{2}}{4}\right)\nonumber \\
 & = & \frac{(-1)^{h}2^{-2h-1}}{h!\Gamma(h+2)}k^{2h+1}\;.\label{eq:J1summedtoprimes-1}\end{eqnarray}
To verify the lowest $43$ summed series for $k\rightarrow8$, one
has to take the upper limit on the number of terms in the series $\geq h+20$
in order to obtain a percent difference between left- and right-hand
sides that was $\leq10^{-33}$ because the first \emph{h} terms in
the series do not contribute. For $k\rightarrow5$, this reduces somewhat
to $\geq h+16$. For $h=0$, one needs 23 terms and 20 terms, respectively. 

In the prior paper \cite{stra24b} we noted that because the modified
Bessel functions of the first kind $I_{N}\left(kx\right)$ are related
to the ordinary Bessel functions by the relation \cite{GR5 p. 987 No. 8.511.4}
(p. 961 No. 8.406.3)

\begin{equation}
I_{n}(z)=i^{-n}J_{n}(iz)\;,\label{eq:I_from_J}\end{equation}
one merely needs to multiply by $i^{-n}$ and set $k=i$ in (\ref{eq:a_L_as_2F3})
to obtain the $I_{0}\left(x\right)$ Fourier-Legendre series. Furthermore,
one sees that $I_{0}$ expressed in powers of x is simply the $J_{0}$
version with all of the negative signs reversed. This is not true
of (\ref{eq:Fourier-Legendre series}) because the arguments of the
Legendre polynomials do not undergo $x\rightarrow ix$ since they
derive from the definition of the Fourier-Legendre series, (\ref{eq:Fourier-Legendre series}).
The \emph{k}-dependence is entirely within the coefficients $a_{LN}\left(k\right)$.

Therefore, the $I_{0}$ Legendre series expansion leads to no new
set of summed series since these would simply be (\ref{eq:1/primes^n=00003D00003D00003D})
with $k=i\kappa$. This is also the case for a Chebyshev expansion.
Clenshaw \cite{Clenshaw} confirms this for $h=0$ on pp. 34-5.

\section{Summed series involving $\,_{1}F_{2}$ hypergeometric functions from
Gegenbauer polynomial expansions of Bessel functions}

One can also use Wimp's Jacobi expansion \cite{jacobi_cheb_Wimp} to
find Gegenbauer polynomial expansions of Bessel functions, since \cite{GR5 p. 987 No. 8.511.4}
(p. 1060 No. 8.962.4)\begin{equation}
P_{2n}^{\left(\lambda-\frac{1}{2},\lambda-\frac{1}{2}\right)}(z)=\frac{\left(\lambda+\frac{1}{2}\right)_{2n}C_{2n}^{\lambda}(z)}{(2\lambda)_{2n}}\;.\label{eq:Jacobi_Gegenbauer}\end{equation}
Unlike those in the first section, the following expansion applies
to non-integer indices as well:

\begin{equation}
J_{\upsilon}\left(kx\right)=\left(kx\right)^{\upsilon}\sum_{L=0}^{\infty}b_{L\upsilon}\left(k\right)C_{2L}^{\lambda}(x)\:,\quad-1\leq x\leq1\label{eq:Gegenbauer series}\end{equation}
where where the coefficients are given by the orthogonality of the
Chebyshev polynomials, \begin{equation}
b_{L\upsilon}\left(k\right)=\frac{2^{2\lambda-1}(2L)!(2L+\lambda)\Gamma(\lambda)^{2}}{\pi\Gamma(2L+2\lambda)}\int_{-1}^{1}\left(J_{\upsilon}\left(kx\right)\left(kx\right)^{-\upsilon}\right)\left(1-x^{2}\right)^{-\frac{1}{2}+\lambda}C_{2L}^{\lambda}(x)\, dx\;,\label{eq:Gegenbauer coefficients}\end{equation}
as

\begin{equation}
\begin{array}{ccc}
b_{L\upsilon}\left(k\right) & ={\displaystyle \frac{(-1)^{L}a^{2L}2^{2L-\nu}\left(\lambda+\frac{1}{2}\right)_{2L}}{\sqrt{\pi}(2\lambda)_{2L}(2L+2\lambda)_{2L}\left(L+\frac{1}{2}\right)_{\nu+\frac{1}{2}}}} & \,_{1}F_{2}\left(L+\frac{1}{2};2L+\lambda+1,L+\nu+1;-\frac{a^{2}}{4}\right)\end{array}\;,\label{eq:b_L_as_1F2 Geg}\end{equation}

The first 22 terms in the Gegenbauer polynomial expansions of $J_{0}\left(kx\right)$
(\ref{eq:Fourier-Legendre series}) are then, with $k=1$ and arbitrarily
taking $\lambda=\frac{1}{4},$ are%
\begin{comment}
from the file a\_ \{L, NN\} from integral and series m7.nb%
\end{comment}
{} \begin{eqnarray}
J_{0}\left(x\right) & \cong & 0.904078771191585521024227636544096\, C_{0}^{\frac{1}{4}}(x)\nonumber \\
 & - & 0.377480902332903752477356198652003\, C_{2}^{\frac{1}{4}}(x)\nonumber \\
 & + & 0.00985645918454006348253321451683292\, C_{04}^{\frac{1}{4}}(x)\nonumber \\
 & - & 0.0000929144245327682841642709978007872\, C_{6}^{\frac{1}{4}}(x)\nonumber \\
 & + & 4.51192238929050409752370668969243\times10^{-7}\; C_{8}^{\frac{1}{4}}(x)\nonumber \\
 & - & 1.33557953986611692879627373122257\times10^{-9}\; C_{10}^{\frac{1}{4}}(x)\nonumber \\
 & + & 2.66185765952711910618049951726347\times10^{-12}\; C_{12}^{\frac{1}{4}}(x)\nonumber \\
 & - & 3.81525698311458688534308130184138\times10^{-15}\; C_{14}^{\frac{1}{4}}(x)\nonumber \\
 & + & 4.12174698882181605290995488668659\times10^{-18}\; C_{16}^{\frac{1}{4}}(x)\nonumber \\
 & - & 3.47649878544013257006577318271996\times10^{-21}\; C_{18}^{\frac{1}{4}}(x)\nonumber \\
 & + & 2.35284611436757064520926642520417\times10^{-24}\; C_{20}^{\frac{1}{4}}(x)\nonumber \\
 & - & 1.30601535036874068380434807654702\times10^{-27}\; C_{22}^{\frac{1}{4}}(x)\nonumber \\
 & + & 6.05330821302322601332159315076677\times10^{-31}\; C_{24}^{\frac{1}{4}}(x)\nonumber \\
 & - & 2.37803900637432785965238868426667\times10^{-34}\; C_{26}^{\frac{1}{4}}(x))\nonumber \\
 & + & 8.01904904818064037541914772609834\times10^{-38}\; C_{28}^{\frac{1}{4}}(x)\nonumber \\
 & - & 2.34648500711595153019299896447757\times10^{-41}\; C_{30}^{\frac{1}{4}}(x)\nonumber \\
 & + & 6.01437661018790357782076353076573\times10^{-45}\; C_{32}^{\frac{1}{4}}(x)\nonumber \\
 & - & 1.36150461807454631533129344677808\times10^{-48}\; C_{34}^{\frac{1}{4}}(x)\nonumber \\
 & + & 2.74196263898788782515776484348033\times10^{-52}\; C_{36}^{\frac{1}{4}}(x)\nonumber \\
 & - & 4.94454559738665023143625856030709\times10^{-56}\; C_{38}^{\frac{1}{4}}(x)\nonumber \\
 & + & 8.03022232133135996468784524426669\times10^{-60}\; C_{40}^{\frac{1}{4}}(x)\nonumber \\
 & - & 1.18066972928855334355199708640780\times10^{-63}\; C_{42}^{\frac{1}{4}}(x)\quad-1\leq x\leq1\;.\label{eq:J0 in Chebyshev-1}\end{eqnarray}

At the upper limit of applicability, $x=1$, this gives 33-digit accuracy,
$J_{0}\left(1\right)=0.765197686557966551449717526102663$. (Even
at $x=8$, this gives a result accurate to 15-digits, $0.171650807137554$.)
One gets a different representation with similar accuracy if one takes
$\lambda=4$. The convergence is not any faster for either than for
the Chebyshev version (\ref{eq:J0 in Chebyshev x/8}), so there is
no strong motivation for programmers to switch to this representation
of Bessel functions from the well-established computer codes for Chebyshev expansions. 

There is, however, interesting research into the utility of Gegenbauer
expansions in an analytical context. To note just three examples,
Bezubik, D\k{a}browska, and Strasburger \cite{Bezubik_et_al} derive
an expansion of plane waves $e^{ir\left(\xi|\eta\right)}$ into an
infinite series over $i^{m}\left(\alpha+m\right)J_{\alpha+m}\left(r\right)C_{m}^{\alpha}(\left(\xi|\eta\right))$,
and Elgindy and Smith-Miles \cite{Elgindy and Smith-Miles} develop
a numerical quadrature based on a truncated Gegenbauer expansion
series. A third example is Jens Keiner's method \cite{Keiner} of converting
from one expansion in $C_{j}^{\alpha}(x)$ to another expansion in
$C_{j}^{b}(x)$.

However, neither the relative numerical utility of Gegenbauer expansions,
nor expansions in an analytical context are the morivations for this
paper. We instead seek to
sum additional infinite series involving $\,_{1}F_{2}$ hypergeometric
functions. To extract the powers, we use the conversion for Gegenbauer
polynomials that is equivalent to (\ref{eq:PBM3 p. 468 No. 7.3.1.207}),
which is  \cite{GR5 p. 987 No. 8.511.4} (GR5 p. 1051 No. 8.932.2)%
\begin{comment}
from the file J0 as GegenbauerC m7.nb%
\end{comment}
{}

\begin{equation}
C_{2L}^{\lambda}(x)=\frac{(-1)^{L}}{(L+\lambda)B(\lambda,L+1)}\,_{2}F_{1}\left(-L,L+\lambda;\frac{1}{2};x^{2}\right)\equiv\frac{(-1)^{L}}{(L+\lambda)B(\lambda,L+1)}\sum_{m=0}^{\infty}\frac{\left(x^{2}\right)^{m}(-L)_{m}(L+\lambda)_{m}}{m!\left(\frac{1}{2}\right)_{m}}\:.\label{eq:GR5 p. 1051 No. 8.932 .2}\end{equation}

The summed series associated with the power $x^{2h}$ may simply be
written as

\begin{eqnarray}
\sum_{L=0}^{\infty}\, &  & \hspace{-2.9cm}\frac{(-2)^{2L}(-L)_{h}\left(\lambda+\frac{1}{2}\right)_{2L}(L+\lambda)_{h}}{\sqrt{\pi}h!\left(\frac{1}{2}\right)_{h}\left(L+\frac{1}{2}\right)_{\frac{1}{2}}(L+\lambda)(2\lambda)_{2L}(2L+2\lambda)_{2L}B(\lambda,L+1)}k^{2L}\,_{1}F_{2}\left(L+\frac{1}{2};L+1,2L+\lambda+1;-\frac{k^{2}}{4}\right){}_{,}\nonumber \\
 & =\frac{(-1)^{h}2^{-2h}k^{2h}}{h!\Gamma(h+1)} & \;.\label{eq:J0_C_series}\end{eqnarray}

This has an identical right-hand side as for the Legendre (\ref{eq:1/primes^n=00003D00003D00003D})
and Chebyshev (\ref{eq:J0_T_series}) versions, and it holds for every
value of $\lambda$. That is we have just summed an infitite set of
infitite sets of infinite series involving $\,_{1}F_{2}$ hypergeometric
functions. To see how this plays out in practive, consider two extreme
values, $\lambda=2^{\pm20}$. For $h=1$ (associated with $x^{2}$) and $\lambda=2^{-20}$, the
first eight terms sum as

\begin{equation}
\begin{array}{c}
0\\
-0.234776027081720679198861338236978 \\
-0.0149856953860168611951004494702182 \\
-0.000236617512932378657466894715711978 \\
-1.653516950294282858347187372908306\times10^{-6} \:.\\
-6.486087810927263603219843086719685\times10^{-9} \\
-1.62636408715893081661576487444697\times10^{-11} \\
-2.82986734360173162046207736379133\times10^{-14} -\cdots\\
=-0.24999999999999996 \end{array}\label{eq:2^-20sum}\end{equation}

For $\lambda=2^{20}$, the second term is almost sufficient by itself:

\begin{equation}
\begin{array}{c}
0\\
-0.249999955296648756645694568085746 \\
-4.470334680250094078684477109090418\times10^{-8} \\
-4.44085305583632491122417076683932\times10^{-14} \\
-3.08808728006191746252089269280760\times10^{-22} \:.\\
-1.65656014384509551797626946530621\times10^{-29} \\
-7.24073285415562941500526940820453\times10^{-37} \\
-2.67164795422661275973040680548870\times10^{-44} -\cdots\\
=0.250000000000000000000000000000000 \end{array}\label{eq:2^-20sum-1}\end{equation}

We turn now to display the first 22 terms in the Gegenbauer polynomial
expansions of $J_{1}\left(kx\right)$ (\ref{eq:Fourier-Legendre series}),
with $k=1$ and arbitrarily taking $\lambda=\frac{1}{4},$ are%
\begin{comment}
from the file a\_ \{L, NN\} from integral and series m7.nb%
\end{comment}
{} \begin{eqnarray}
J_{1}\left(x\right) & \cong & 0.475683429275416807386224265471041\; x\, C_{0}^{\frac{1}{4}}(x)\nonumber \\
 & - & 0.0962237678006581825132637018597388\; x\, C_{2}^{\frac{1}{4}}(x)\nonumber \\
 & + & 0.00165923280553475766418121861007493\; x\, C_{04}^{\frac{1}{4}}(x)\nonumber \\
 & - & 0.0000116849150281699572948996291216163\; x\, C_{6}^{\frac{1}{4}}(x)\nonumber \\
 & + & 4.5303088506394388853845501270703\times10^{-8}x\; C_{8}^{\frac{1}{4}}(x)\nonumber \\
 & - & 1.11623410748844105451625882776928\times10^{-10}x\; C_{10}^{\frac{1}{4}}(x)\nonumber \\
 & + & 1.90550296957009549791418728733899\times10^{-13}x\; C_{12}^{\frac{1}{4}}(x)\nonumber \\
 & - & 2.38861692204435794092836019335553\times10^{-16}x\; C_{14}^{\frac{1}{4}}(x)\nonumber \\
 & + & 2.29299953783708159991903279787185\times10^{-19}x\; C_{16}^{\frac{1}{4}}(x)\nonumber \\
 & - & 1.74020202094491142079186625047494\times10^{-22}x\; C_{18}^{\frac{1}{4}}(x)\nonumber \\
 & + & 1.07047764587989141691542634270970\times10^{-25}x\; C_{20}^{\frac{1}{4}}(x)\nonumber \\
 & - & 5.44604885209780265146726077614161\times10^{-29}x\; C_{22}^{\frac{1}{4}}(x)\nonumber \\
 & + & 2.32978017343783698464445765163641\times10^{-32}x\; C_{24}^{\frac{1}{4}}(x)\nonumber \\
 & - & 8.49800957174229357388497217989335\times10^{-36}x\; C_{26}^{\frac{1}{4}}(x))\nonumber \\
 & + & 2.67439765035844790866837011922870\times10^{-39}x\; C_{28}^{\frac{1}{4}}(x)\nonumber \\
 & - & 7.33611068017397602824622074628328\times10^{-43}x\; C_{30}^{\frac{1}{4}}(x)\nonumber \\
 & + & 1.76965188025225356497750305500631\times10^{-46}x\; C_{32}^{\frac{1}{4}}(x)\nonumber \\
 & - & 3.78333047168286059388389868568285\times10^{-50}x\; C_{34}^{\frac{1}{4}}(x)\nonumber \\
 & + & 7.21804990626747371147788669564168\times10^{-54}x\; C_{36}^{\frac{1}{4}}(x)\nonumber \\
 & - & 1.23650210830378663827086788837057\times10^{-57}x\; C_{38}^{\frac{1}{4}}(x)\nonumber \\
 & + & 1.91247206300887832512635377246951\times10^{-61}x\; C_{40}^{\frac{1}{4}}(x)\nonumber \\
 & - & 2.68399957497958828507307548313041\times10^{-65}x\; C_{42}^{\frac{1}{4}}(x)\quad-1\leq x\leq1\;.\label{eq:J0 in Chebyshev-1-1}\end{eqnarray}

At the upper limit of applicability, $x=1$, this gives 33-digit accuracy,
$J_{1}\left(1\right)=0.440050585744933515959682203718915$. (Even
at $x=8$, this gives a result accurate to 16-digits, $J_{1}\left(8\right)=0.2346363468539146$.)
One gets a different representation with similar accuracy if one takes
$\lambda=4$. The convergence is not any faster for either than for
the Chebyshev version (\ref{eq:J0 in Chebyshev x/8}), so there is
no strong motivation for programmers to switch to this representation
of Bessel functions. 

The summed series derived from Gegenbauer polynomial expansions of
$J_{\upsilon}\left(x\right)$, may be found for any value of $\nu$,
not just integer values, given that it is derived from Wimp's Jacobi
expansion \cite{jacobi_cheb_Wimp}. The series associated with the
power $x^{2h+\nu}$ in the general-$\nu$ case may be written as

\begin{eqnarray}
\sum_{L=0}^{\infty}\, &  &  \hspace{-3.5cm} \frac{(-1)^{2L}2^{2L-\nu}(-L)_{h}\left(\lambda+\frac{1}{2}\right)_{2L}(L+\lambda)_{h}}{\sqrt{\pi}h!\left(\frac{1}{2}\right)_{h}(L+\lambda)(2\lambda)_{2L}(2L+2\lambda)_{2L}\left(L+\frac{1}{2}\right)_{\nu+\frac{1}{2}}B(\lambda,L+1)}k^{2L+\nu}\,_{1}F_{2}\left(L+\frac{1}{2};2L+\lambda+1,L+\nu+1;-\frac{k^{2}}{4}\right){}_{,}\nonumber \\
 & =\frac{(-1)^{h}2^{-2h-\nu}k^{2h+\nu}}{h!\Gamma(h+\nu+1)} & \;.\label{eq:J0_C_nu_series}\end{eqnarray}
Thus, we have just summed an infitite set of infitite sets of doubly
infinite series involving $\,_{1}F_{2}$ hypergeometric functions
since the expression holds for every value of $\lambda$ and holds
for every value of $\nu$.

Comparison of (\ref{eq:J0_C_series}) and (\ref{eq:J0_C_nu_series})
shows us how to extend the Chebyshev (\ref{eq:J0_T_series}) version
to every value of $\nu$: 

\begin{eqnarray}
\sum_{L=1}^{\infty}\, &  & \hspace{-0.9cm}\frac{(-1)^{L}2^{-2L-\nu}\left(\frac{1}{2}-L\right)_{L-h}(-L)_{L-h}}{L!(L-h)!(1-2L)_{L-h}\Gamma(L+\nu+1)}k^{2L+\nu}\,_{1}F_{2}\left(L+\frac{1}{2};2L+1,L+\nu+1;-\frac{k^{2}}{4}\right)\nonumber \\
 & = & \frac{(-1)^{h}2^{-2h-1}}{h!\Gamma(h+2)}k^{2h+1}\;.\label{eq:J_nu_T_series}\end{eqnarray}
This is a modest variation on the form one gets by simply setting
$\lambda=0$ in (\ref{eq:J0_C_nu_series}), since $T_{\nu}(z)=\frac{1}{2}\nu C_{\nu}^{0}(z)$.

An extension of the Legendre sets (\ref{eq:1/primes^n=00003D00003D00003D})
and (\ref{eq:J1summedtoprimes}) to larger integer values of $\nu$
is not obvious, but one can get such a form directly from (\ref{eq:J0_C_nu_series})
for all values of $\nu$ since $P_{\nu}(z)=C_{\nu}^{\frac{1}{2}}(z)$.

\section{Conclusion}

I have shown how to sum doubly infinite sets of infinite series involving
$\,_{1}F_{2}$ hypergeometric functions, derived from Chebyshev polynomial
expansions of Bessel functions, and trebly infinite sets of infinite
series involving $\,_{1}F_{2}$ hypergeometric functions from Gegenbauer
polynomial expansions of Bessel functions of the first kind $J_{\nu}\left(kx\right)$.
The utility of any one of these summed series for future researchers
is, of course, not guaranteed, but given the relative paucity of infinite
series whose values are known (e.g., 24 pages in Gradshteyn and Ryzhik
compared to their 900 pages of known integrals), one hopes that adding
a trebly infinite set of infinite series of $\,_{1}F_{2}$ functions
whose values are now known will be of use to some. 

\vspace{0.2cm}
 \textbf{Funding:} This research received no external funding.

\vspace{0.2cm}
 \textbf{Data Availability:} Data are contained within the article.
 
\vspace{0.2cm}
 \textbf{Conflicts of Interest:} The author declares no conflicts
of interest.

\end{document}